\documentclass[11pt]{amsart}
\usepackage{geometry}                
\usepackage{graphicx}
\usepackage{amssymb}
\usepackage{epstopdf}
\usepackage{tikz}
\newtheorem{Thm}{Theorem}

\newtheorem{Def}{Definition}
\newtheorem{Cor}{Corollary}

\title[Arrow's Theorem and the TSP]{Connecting Arrow's Theorem, Voting Theory, and the Traveling Salesperson Problem}
\author[D.~Saari]{Donald G. Saari}  
\thanks{This paper is a written version of my 4/8/2022 JMM presentation; these results provide applications of  conclusions  in \cite{tandd, new}.   My thanks to  George Hazelrigg for our several  discussions of these topics.  This work is part of a National Science Foundation project under NSF Award Number CMMI-1923164.}

\address{IMBS; University of California, Irvine\\ Irvine, California 92617-5100}
\email{dsaari@uci.edu}

\begin{document}
\maketitle

\begin{abstract}  Problems with majority voting over pairs as represented by Arrow's Theorem   and those of finding the lengths of  closed paths as captured by  the  Traveling Salesperson Problem (TSP) appear to have nothing in common.  In fact, they are connected.  As shown,  pairwise voting and a version of the TSP     share  the same domain where each  
system can be simplified  
by   restricting it  to
  complementary  regions to  eliminate  extraneous  terms.  Central for doing so is the Borda Count, where it is shown that its outcome most accurately reflects the  voter preferences. 
\end{abstract}

\section{Introduction}

Among the many challenges  posed by discrete mathematics and the Social Sciences  are  aspects of voting theory as characterized by Arrow's Theorem and the properties of closed paths on a graph as typified  by  the Traveling Salesperson Problem (TSP).  Surprisingly,   
  both topics   can be analyzed with  essentially the same  
  approach.  After introducing the commonality with  
    examples,  
it is shown how   each situation can be simplified     by emphasizing   different regions of the associated  geometry. 

\subsection{Voting theory and Arrow's Theorem} Central to  Arrow's Theorem \cite{arrow} is his IIA condition, which  requires  a profile's  conclusion to be completely  determined by the rankings of  
 its associated  paired comparisons.\footnote{Even stronger,  a decision method  satisfies ``Independence of irrelevant alternatives," IIA, iff it  can be expressed in terms of independent paired comparison methods \cite{bk, pc}.  This means it is not a negative feature if a decision  method, such as the plurality vote,   fails to satisfy IIA; it  only means that the method  cannot be so expressed. Conversely, any method that cannot be expressed in this manner does not satisfy IIA.  Thus, Arrow's Theorem is a result about the problems that accompany  paired comparisons.}   A $N=70$ person majority vote example over    the  $n=3$ alternatives  $\{A_j\}_{j=1}^3$  is where  
 \begin{equation}\label{ex: 1}  25  \textrm{ prefer }A_1\succ A_2\succ A_3, \,  23 \textrm{ prefer } A_2\succ A_3\succ A_1, \,  22 \textrm{ prefer }A_3 \succ A_1\succ A_2. \end{equation}  This leads to the cyclic pairwise majority vote outcomes of  \begin{equation}
 \label{eq: introductory}
A_1\succ A_2 \textrm{ by 47:23, } A_2\succ A_3 \textrm{ by 48:22, and } A_3\succ A_1 \textrm{ by 45:25},\end{equation}  which violate  Arrow's objective of obtaining  a transitive conclusion.  This difficulty reflects  
the long-standing objective  in voting theory, which is to replace  cyclic  outcomes  with appropriate transitive rankings,  or   at least with outcomes that identify  a  ``best choice."    Prominent  approaches  were developed by the mathematicians Dodgson \cite{dodgson} in 1876 and    Kemeny \cite{kemeny} in 1959.

A way to compare  how    an alternative fares in a paired comparison   is to compute how its tally differs from the average score of $ \frac{\textrm{number of voters}}2 =\frac N2$.  (In what follows, $N$ is the number of voters, $n$ is the number of alternatives.)  So  
\begin{equation}\label{def: dij} d_{i, j} = \frac12[ A_i\textrm{'s tally} -  A_j\textrm{'s tally}] =A_i\textrm{'s tally} -  \frac N2  = -d_{j, i}.\end{equation}  
For the Eq.~\ref{eq: introductory}  values, where $N=70$ and $n=3,$ \begin{equation}\label{eq: dijvotediff} d_{1, 2} = 47-\frac{70}2 = 12=-d_{2, 1}, \,  d_{2, 3} = 13=-d_{3, 2}, \textrm{ and } d_{1, 3} = -10 = -d_{3, 1}.\end{equation} 

\subsection{TSP}  To introduce the notation for the 
    TSP system  discussed here,  it takes 30 minutes to walk from home, H,  to campus, C,   while returning uphill takes 40.  As the average  is 35 minutes,  
 returning home requires  5   minutes above the average    denoted by $C  \stackrel{5}{\longrightarrow}  H$;   traveling   in the opposite direction takes 5 fewer minutes or    $H \stackrel{-5}{\longrightarrow} C$.  More generally, if $d_{i, j}$ represents the   ``difference from the average"      cost of going  from $A_i$ to $A_j$, then  \begin{equation}\label{eq: path dij}
  A_i \stackrel{d_{i, j}}{\longrightarrow}  A_j \textrm{ and }  A_j \stackrel{-d_{i, j}}{\longrightarrow}  A_i \textrm{ are equivalent},\end{equation} 
  or, as true with Eq.~\ref{def: dij},  \begin{equation}\label{eq: dijforcost}  d_{i, j} = -d_{j, i}.\end{equation} 
  
 
\begin{tikzpicture}[xscale=0.5, yscale=0.52]

\draw[<-] (1.18, 0) -- (3.9, 0); 
\node[below] at (1, 0) {$A_1$};
\node[below] at (2.5, 0) {3};
\draw[->] (1.15, .06) -- (4.9, 2);   
\node at (1.9, .5) {\small 3};
\draw[->] (1.15, .15) -- (3.9, 3.9); 
\node at (2.1, 1.5)  {\small 2};
\draw[<-] (1, .2) -- (1, 3.9);  
\node at  (1, 3) {\small 2}; 
\node[left] at (0, 2) {$A_6$}; 
\draw[->] (.9, .1) -- (.1, 1.9);  
\node[left] at (.5, 1) {3};
\draw[->] (4.1, .1) -- (4.9, 1.9); 
\node[right] at  (4.5, 1) {11}; 
\node[below] at (4, 0) {$A_2$};
\draw[->] (4, .1) -- (4, 3.85);  
\node at (4, 1) {\small 11};
\draw[dashed]  (3.9, .1) -- (1.1, 3.9); 
\node at (2.1, 2.5) {\small 0};
\draw[->] (3.82, .07) -- (.23, 1.95);  
\node at (.56, 1.7) {\small 8};
\draw[<-] (4.9, 2.1) -- (4.1, 3.9);  
\node[right] at (4.5, 3) {2}; 
\node[right] at (5, 2){$A_3$};
\draw[<-] (4.7, 2.1) -- (1.24, 3.86); 
\node at (4.3, 2.35) {\small 14};
\draw[<-] (4.7, 2) -- (.2, 2); 
\node at (3.3, 2) {\small 3};
\draw[<-] (3.8, 4) -- (1.2, 4); 
\node[above] at (2.5, 4) {9};
\node[above] at (4, 4) {$A_4$}; 
\draw[<-] (3.75, 3.85) -- (.2, 2.1);  
\node at (3.05, 3.4)   {\small 1};
\node[above] at (1, 4) {$A_5$}; 
\draw[->]  (.83, 3.9) -- (.1, 2.1); 
\node[left] at (.5, 3) {8};

\node[below] at (2, -1) {{\bf Figure 1.} $\mathcal G_A^6$};

\end{tikzpicture}

Figure~1 is a  typical TSP example, where the  graph catalogues all   ``differences from averages" information about the  alternatives $\{A_j\}_{j=1}^6$ and its ${6\choose2} = 15$ pairs.  (Subscripts ``$A$" and `$S$"  
refer, respectively,   to whether the graph has an asymmetric or a  symmetric structure.)   Only positive cost directions are  
 displayed because (Eqs.~\ref{eq: path dij}, \ref{eq: dijforcost}) traveling counter to an arrow is a negative, or  below average cost.  A standard  TSP objective  is to discover the longest and shortest Hamiltonian circuits.  To appreciate what will be developed, before reading more, let me ask the reader to find the longest such path in Fig.~1.  Recall,  a Hamiltonian circuit is a closed path that starts and ends at a selected vertex $A_j$ and passes  through each of the other vertices  once.  
To illustrate,  the Fig.~1 Hamiltonian path  $$ A_1  \stackrel{-3}{\longrightarrow} A_2 \stackrel{11}{\longrightarrow} A_4 \stackrel{-1}{\longrightarrow} A_6 \stackrel{3}{\longrightarrow} A_3 \stackrel{-14}{\longrightarrow}  A_5 \stackrel{2}{\longrightarrow} A_1$$ has length $-2$, and its  
 reversal has length $2.$ In fact, thanks to  
 Eqs.~\ref{eq: path dij}, 
 \ref{eq: dijforcost},   the tasks of finding the longest and shortest Hamiltonian paths coincide; the reversal of one is the other.

\subsection{A common space}   The  voting and   TSP  questions    involve the $n\choose2$ pairs defined by  $n$ alternatives where costs/differences between vertices of a $\{i, j\}$ pair are measured by $d_{i, j}=-d_{j,i}$ values.\footnote{ 
Rather than  ``difference from average,"  $d_{i, j}$ can be anything; e.g., a natural choice is  the difference between $A_i$ and $A_j$ values.}  Thus,   both settings share the    $\mathbb R_A^{n\choose2}$ domain,  where $\mathbf d_A^n\in\mathbb R_A^{n\choose2}$ has the form \begin{equation}\label{eq: d}
 \mathbf d_A^n=  (d_{1, 2}, d_{1, 3}, \dots, d_{1, n}; d_{2, 3}, \dots, d_{2, n}; d_{3, 4} \dots; d_{n-1, n});  \, d_{i, j} = -d_{j, i}.\end{equation}  Semicolons indicate changes in the first subindex. The  $\mathbb R_A^{n\choose2}$ structure is developed next.

To start, recall that if a triplet $\{A_i, A_j, A_k\} $  defines a transitive ranking, it can be expressed  with these $d_{s, u}$ values; e.g.,   should $d_{i, j}>0$ and $d_{j, k}>0$, then transitivity requires  
 that  $d_{i, k}>0$.  These inequalities are borrowed from the structure of points on a line where $p_i>p_j$ and $p_j>p_k$ require $p_i>p_k$.  But these points also satisfy the stronger algebraic  relationship $(p_i-p_j)+(p_j-p_k) = (p_i-p_k)$.  The following mimics this equality.   
\begin{Def}\label{def: transitive} \cite{tandd}   Vector $\mathbf d_A^n\in \mathbb R_A^{n\choose2}$  is strongly transitive iff  each triplet $\{i, j, k\}$ satisfies \begin{equation}\label{def: st} d_{i, j}+d_{j,k} = d_{i, k}.\end{equation}.  \end{Def}  

This condition was introduced for profiles in \cite{et-1}; the  decision theory version used here  comes from \cite{tandd}.  The subspace of strongly transitive vectors, $\mathbb{ST}_A^n$,  is described next.

\begin{Thm}\label{thm: st}  \cite{tandd} The set of strongly transitive $\mathbf d_{A, st}^n \in \mathbb R_A^{n\choose2}$, denoted by $\mathbb {ST}^n$,\footnote{In \cite{et-1, kem-geometry}, $\mathbb{ST}_A^n$ is called the ``transitivity plane."} is a $(n-1)$-dimensional linear subspace of $\mathbb R_A^{n\choose2}$.\end{Thm}

While details are in \cite{tandd}, proving that $\mathbb {ST}_A^n$ is a linear subspace is a standard exercise.  The assertion about its dimension
 follows from the fact that   
  any $d_{j,k}$ can be expressed as
$d_{j, k} = d_{j, 1} + d_{1, k}$ (Eq.~\ref{def: st}),  
 so  all $d_{j, k}$ values for $\mathbf d^n_{A, st}\in \mathbb{ST}^n_A$ can be determined  from the $(n-1)$ terms $\{d_{1, s}\}_{s=2}^n.$  $\square$\smallskip

 The dimensions of $\mathbb R_A^{n\choose2}$ and $\mathbb{ST}_A^n$ dictate that $\mathbb C_A^n$, the orthogonal complement  of $\mathbb{ST}_A^n$,  has dimension ${n-1}\choose2$.  To motivate  the form of these   orthogonal vectors,   expressing Eq.~\ref{def: st} as $x+y=z$, or $x+y-z=0$, identifies  $(1, 1, -1)$ as a normal vector. Thus the $x=y=-z=1$ values define a normal  vector $d_{i, j}=  d_{j, k}= d_{k, i} =1$, which is  represented by the cyclic   
 $A_i\succ A_j, \, A_j\succ A_k, \, A_k\succ A_i$  
  where  the  differences between values is the same constant.

\begin{Thm} \label{thm: cyclic}  \cite{tandd} The ${n-1}\choose2$ dimensional linear subspace $\mathbb C_A^n$ consists of cycles where one   basis, which consists  of three-cycles,  is
 $\{\mathbf d^n_{1, j, k}\}_{1<j<k\le n}$.  The only non-zero terms of $\mathbf d^n_{1, j, k}$  are $d_{1, j} = d_{j, k} = d_{k, 1}=1.$\end{Thm}
 
 According to these  theorems, $\mathbb R_A^{n\choose2}$ nicely separates into  linear subspaces of strongly transitive vectors,  $\mathbb{ST}_A^n$,  and  
  cyclic vectors, $\mathbb C_A^n$, where $\mathbb C_A^n$ has a basis consisting of a special type of three-cycles. This means that $\mathbf d^n_A \in \mathbb R_A^{n\choose2}$ has a unique decomposition \begin{equation}\label{eq: forgotten} \mathbf d^n_A = \mathbf d^n_{A, st} + \mathbf d^n_{A, cyclic} \end{equation}
   where $\mathbf d^n_{A, st}$ and $\mathbf d^n_{A, cyclic}$ are, respectively, the orthogonal projection of $\mathbf d^n_A$ to $\mathbb{ST}^n$ and $\mathbb C^n_A$.   
 
  An immediate consequence of Eq.~\ref{eq: forgotten} is that  
  a transitive vector, which    is not strongly transitive, is a hybrid.  
 
 \begin{Cor}\label{cor: transitive} For $n\ge 3$, if $\mathbf d_A^n$ is transitive but not strongly transitive, then there are unique non-zero  vectors  $\mathbf d^n_{A, st}\in \mathbb {ST}_A^n$ and $\mathbf d^n_{A, cyclic}\in \mathbb C_A^n$ so that $\mathbf d_A^n=\mathbf d^n_{A, st} + \mathbf d^n_{A, cyclic}$. \end{Cor} 
 
 What limits the use of 
transitive vectors is that  they fail to define a linear subspace.   For instance,
   $d_{1, 2} = d_{2, 3}=2$, $d_{1, 3}=1$  
    and  $d^{*}_{1, 2} = d^{*}_{2, 3}=-1$, $d^{*} _{1, 3}=-2$ are transitive triplets, but their sums (e.g., $d_{i, j} + d^*_{i, j} = d^+_{i, j})$ define the cyclic $ d^+_{1, 2} = d^+_{2, 3} = 1$, $d^+_{1, 3} = -1$.   
According to Cor.~\ref{cor: transitive}, it is reasonable to treat a transitive vector  
 as  a strongly transitive choice that is lightly  contaminated with cyclic terms.  The following theorem partly explains the ``lightly contaminated" modifier by showing that rankings of  transitive $\mathbf d_A^n$ and its strongly transitive component $\mathbf d^n_{A, st}$ can differ, but not   radically.   This holds even for a non-transitive  $\mathbf d_A^n$ that   has,  at least, a Condorcet winner and  loser.  (A Condorcet winner  is a candidate who beats each of the other candidates  in  majority vote comparisons.  A Condorcet loser loses all pairwise majority votes. They can exist even without transitivity; e.g.,   
  for $n=5,$ $A_1$ could be the Condorcet winner, $A_5$ the Condorcet loser, and $A_2, A_3, A_4$ define a cycle.)

\begin{Thm}\label{thm:  ranking} If $\mathbf d_A^n$ has a unique Condorcet winner $A_1$ and a unique Condorcet loser $A_n$, then $A_1$ is strictly ranked above  $A_n$ in $\mathbf d^n_{A, st}$.\end{Thm} 

The converse  is to determine what happens by adding cyclic terms  to a $\mathbf d^n_{A, st}.$  If the resulting  $\mathbf d_A^n$ is wildly cyclic, not much can be stated.  But if  $A_1$ and $A_n$ are, respectively, the top and bottom ranked candidates of the  $\mathbf d^n_{A, st}$ and if  $\mathbf d_A^n$ remains transitive, the question is whether  
 $A_1$ must be ranked above $A_n$ in $\mathbf d_A^n=\mathbf d^n_{A, st}+\mathbf d^n_{A, cyclic}$.    Proofs of this kind of  results  can be messy.  But  to illustrate  Thm.~\ref{thm: cyclic}, the basic idea is developed in Sect.~\ref{sect: proof}  for $n=3, 4$.
The details of these proofs are similar to   
 earlier relationships that were developed about the Kemeny and paired voting rankings relative to the Borda ranking \cite{et-1, kem-geometry}.

 \section{  Voting Theory}

 As a central objective in voting and decision theory is to obtain transitive outcomes,   
   the cyclic components of  $\mathbf d_A^n$    introduce problems.   
A resolution is obvious; drop the troubling $\mathbf d^n_{A, cyclic}$ cyclic term and retain  only  $\mathbf d^n_{A, st}$.    Dodgson's method partly  does so  
  by replacing  $\mathbf d_A^n$   
 with a  vector that may have cyclic terms, but at least  it  has a Condorcet winner.  ($A_j$ is a Condorcet winner iff $d_{j, k}>0$ for all $k\ne j$. That is, $A_j$ is ``better than" all other alternatives in paired comparisons.)     Dodgson's  method \cite{dodgson}, then,   projects  $\mathbf d_A^n\in \mathbb R_A^{n\choose2}$  to   
 the nearest $\mathbb R_A^{n\choose2}$   subset  where all  vectors have a Condorcet winner.  For $n\ge 4$,     these regions  have cyclic components;  
 e.g., for $n=4$, the dimension of $\mathbb{ST}_A^4$ is only three, while the  Condorcet  subset   has the full ${4\choose2}=6$ dimension.  For instance,  
  $A_1$ is the Condorcet winner  in the eight  $\mathbb R_A^{6}$ orthants  where $d_{1, 2}>0, d_{1, 3}>0, d_{1, 4}>0$.  Without  imposing  restrictions on the signs of $d_{2, 3}, d_{2, 4}, d_{3, 4}$,  cyclic behavior is allowed among $\{A_2, A_3, A_4\}$.

 Kemeny  adopted the  
 more ambitious goal \cite{kemeny} of replacing  $\mathbf d_A^n$ with a transitive outcome.  Similar to Dodgson,  he created a projection mapping that sends  $\mathbf d_A^n$  into the nearest  $\mathbb R_A^{n\choose2}$ subset consisting of transitive outcomes.  
 
Both methods provide electoral  relief at a first level,  but  a deeper investigation  reveals  
 a host of  other  subtle difficulties that  
 cast doubt on the reliability of  these approaches. Important results in this direction were found in a series of papers by Ratliff \cite{tommy1, tommy2, tommy3}.   
 For instance,  the above ``projection" descriptions makes it reasonable to expect  
   that the  Dodgson and Kemeny outcomes are related;  perhaps  the Dodgson winner always is the top-ranked Kemeny candidate.  But Ratliff proved that, in general,  such assertions  are false.  A small sample of his findings follows.

\begin{Thm}\label{thm: tommy} \cite{tommy1, tommy2, tommy3} For $n\ge 4$ candidates,  select an integer k, where $1\le k \le n.$  There exist paired comparison  examples where the Dodgson winner is the $k^{th}$ ranked Kemeny winner.  

 Dodgson's method can be generalized to select a committee of $k\ge 2$ candidates by using Dodgson's projection method to  a subset of $\mathbb R_A^{n\choose 2}$ where all vectors have $k$  candidates where each is  ranked above the $(n-k)$ remaining candidates.   For  integers $k$ and $s$ satisfying $1\le s  < n-k,$ $k\ne s$,  
there exist examples where with the Dodgson's generalized approach, no candidate in the Dodgson committee of $s$   is in the   Dodgson committee of  $k$   candidates.\end{Thm}

As Ratliff proved, rather than being top-ranked,  the Dodgson winner can  end up being anywhere in a Kemeny outcome; it can even be bottom-ranked.  Moreover,  the Dodgson winner need not be in a Dodgson-Ratliff  committee of two or three.    These unexpected  conclusions are consequences of   
the  cyclic terms that    remain even  after the Dodgson and Kemeny projections.  For instance, unless Kemeny's outcome  is strongly transitive, it contains cyclic components (Cor.~\ref{cor: transitive}) that can create   other difficulties.   
 
Stated differently,  
  the Dodgson and Kemeny procedures    remove  only as many of the $\mathbf d_A^n$ cyclic components as needed to attain their objectives.   Without a thorough  cleansing,   it is reasonable to  
 expect  other mysterious properties: this happens.   
This also holds for  Arrow's Theorem;   
 its negative conclusion is strictly a consequence of the $\mathbf d_{A, cyclic}^n$  component.  Completely removing $\mathbf d_{A, cyclic}^n$ converts Arrow's assertion  into a positive result \cite{bk, pc}.  Indeed,  
  the ultimate goal for    decision and voting problems  should be to eliminate {\em all  cyclic components of a given $\mathbf d_A^n$}.       Doing so is 
a common mathematical computation.   
 
\subsection{Orthogonal Projection}   The    standard  way to
  eliminate the unwanted   cyclic terms from a given $\mathbf d_A^n$  is with  the  orthogonal projection \begin{equation}\label{def: orthogonal} P:\mathbb R_A^{n\choose2} \to \mathbb {ST}_A^n.\end{equation}   The computations require  
 finding a basis for $\mathbb {ST}_A^n$ (using Eq.~\ref{def: st}) and  carrying out the associated    vector analysis.   The resulting approach  
 follows; see \cite{tandd} for   details and examples.   

 \begin{Def}\label{def: S} For  alternative $A_j$,      
  the weighted   sum is  \begin{equation}\label{eq: def-borda} \mathcal S_A(A_j)=  \sum_{k=1}^n d_{j, k}; \, j=1, \dots, n.  \end{equation}  \end{Def}

With Eq.~\ref{eq: dijvotediff} values, \begin{equation}\label{eq: S values}  S_A(A_1) =  12-10=2, \, S_A(A_2) = -12 +13 =1, \, S_A(A_3) = -13 + 10 = -3.\end{equation}

 \begin{Thm}\label{thm: computing cpi} \cite{tandd} For $\mathbf d_A^n\in \mathbb R_A^{n\choose2}$, the $d_{i, j}$ value in $\mathbf d^n_{A, st}$ is \begin{equation}\label{eq: dij-dn} d_{i, j} = \frac1n [\mathcal S_A(A_i) - \mathcal S_A(A_j)].\end{equation}  
  \end{Thm}\smallskip
  
Of surprise,  this projection is equivalent to the well known Borda Count.    This Borda procedure tallies a $n$-candidate ballot  by assigning $n-j$ points to the $j^{th}$ ranked candidate.  A candidate's Borda tally, $\mathcal B(A_j)$,  is the sum of points assigned  to $A_j$ over all $N$ ballots.  Using the Eq.~\ref{ex: 1} example, $\mathcal B(A_1) = 25(3-1) + 23(3-3) + 22(3-2) = 72$,  $\mathcal B(A_2) = 25(3-2) + 23(3-1) + 22(3-3) = 71$, and $\mathcal B(A_3) = 25(3-3) + 23(3-2) + 22(3-1) = 67.$

\begin{Thm}\label{thm: BC}\cite{tandd} For any $n\ge 3$, the   Eq.~\ref{def: orthogonal} orthogonal projection of $\mathbf d_A^n\in \mathbb R_A^{n\choose2}$ to $\mathbb{ST}_A^n$, which is $\mathbf d^n_{A, st}$,     is equivalent to the Borda Count.  More precisely,   \begin{equation} \mathcal B(A_j) =   (n-1)\frac N2 +  S_A(A_j),  \, j=1, \dots, n.\end{equation}  and \begin{equation}\label{eq: StoB} \mathcal B(A_j) - \mathcal B(A_k) = S_A(A_j) - S_A(A_k) = nd_{j, k}.\end{equation}  \end{Thm} 

\noindent{\em Proof:} As known (e.g., \cite{geometry, et-1, bk}),  a way to compute  $\mathcal B(A_j)$ is to sum   $A_j$'s tallies over  each of its $(n-1)$  majority vote paired comparisons.  (As an example, using the Eq.~\ref{eq: introductory} tallies, in the $\{A_1, A_2\}$ and $\{A_1, A_3\}$ elections, $A_1$ receives, respectively, 47 and 25 votes; this $47+25=72$ value agrees with $A_1$'s above computed  $\mathcal B(A_1) =  72$.)   According to Eq.~\ref{def: dij}, \begin{equation}  
\mathcal B(A_j) = \sum_{k=1,  k\ne j}^n [d_{j, k} + \frac N2]  
 = (n-1)\frac{N}2 +S_A(A_j).\end{equation}   Thus, for voting theory,   $S_A(A_j)-S_A(A_k) = \mathcal B(A_j)-\mathcal B(A_k).$  
  $\square$    \smallskip 
 
Theorem \ref{thm: BC} has an interesting consequence.  To set the stage, consider all of those voting methods where the outcome is determined by assigning a score to each candidate; e.g., this includes almost all standard methods such as all  positional methods, cumulative voting, Approval Voting, etc.  As these scores satisfy strong transitivity,  the outcome is in $\mathbb {ST}^n_A$.   A natural objective is to have an outcome that most accurately reflects the views (i.e., preferences) of the voters.  That is, find the $\mathbb {ST}^n_A$  outcome that is closest to the data, which is the orthogonal projection of $\mathbf d_A^n$.  This is  equivalent to the Borda Count (Thm.~\ref{thm: BC}).   

Combining comments from the above provide new explanations for earlier results. For instance, the following, which has the spirit of Thm.~\ref{thm:  ranking},  was proved by using these tools.

\begin{Thm}\label{thm: kemtoborda} \cite{et-1, kem-geometry} For $n\ge 3$, Kemeny's method ranks a Borda winner over the Borda loser.  Conversely, the Borda Count ranks a Kemeny winner over the Kemeny bottom ranked candidate.

For paired comparison majority votes, the Borda outcome ranks the Condorcet winner over the Condorcet bottom ranked candidate.  Conversely,  is the paired comparisons define a transitive ranking,  the Borda winner is ranked over the Borda bottom ranked candidate.  \end{Thm} 

The Kemeny outcome is a transitive vector, while the strongly transitive  Borda ranking comes from the orthogonal projection of $\mathbf d_A^n$ to $\mathbb {ST}_A^n$.  As  Cor.~\ref{cor: transitive} asserts, a transitive ranking is a strongly transitive ranking   clouded by cyclic terms, which captures the flavor of the Kemeny outcome.

Theorem \ref{thm: BC} provides an explanation why the Borda Count has so many positive properties; examples with decision methods  are in \cite{tandd}.  The following completes the introductory comments about Arrow's Theorem.
\begin{Thm} \cite{geometry, et-1}  If $\mathbf d_A^n\in \mathbb{ST}_A^n$, then the Borda Count and majority vote paired comparisons both satisfy Arrow's Theorem.  The Borda Count is the only positional voting method that satisfies Arrow's conditions.
\end{Thm}   

The last assertion holds because  the Borda Count is the only positional method where  its outcomes are determined by the outcomes of  majority votes over pairs. The outcome for all other positional methods  need not be related, in any manner, to the paired comparison outcomes \cite{dictionary}.    

 \section{Turning to the TSP}
 
 According to  Eq.~\ref{eq: forgotten},   $\mathbf d_A^n\in \mathbb R^{n\choose2}$ can be uniquely expressed as 
 $\mathbf d_A^n=\mathbf d^n_{A, st} + \mathbf d^n_{A, cyclic}$, where $\mathbf d^n_{A, st} \in \mathbb {ST}_A^n$ and $\mathbf d^n_{A, cyclic}\in \mathbb C_A^n$. 
 As described above, voting and decision theories seek transitive outcomes, which means that the cyclic  $\mathbf d^n_{A, cyclic}$ component   imposes obstacles.   
  The natural resolution   
 is to eliminate  $\mathbf d^n_{A, cyclic}$  and retain $\mathbf d^n_{A, st}$  by projecting   $\mathbf d_A^n$ into $\mathbb {ST}_A^n.$  
 
 In contrast,   TSP and other closed path concerns  involve {\em  cycles}, so the  linear behavior  
 now is what  creates barriers.  This change in  the objective    converts   
  the strongly transitive $\mathbf d^n_{A, st}$  from being  the desired component into  
  a troublemaker.  Thus,     to analyze    TSP issues,  mimic what was done for voting by     orthogonally  projecting $\mathbf d_A^n$ into $\mathbb C_A^n$ to eliminate the $\mathbf d^n_{A, st}$ term and  retain $\mathbf d^n_{A, cyclic}$.

 What simplifies  finding  $\mathbf d^n_{A, cyclic}$   is that $\mathbf d^n_{A, cyclic}= \mathbf d_A^n-\mathbf d^n_{A, st}$, and $\mathbf d^n_{A, st}$,  which is equivalent to the Borda Count (Thm.~\ref{thm: BC}), is easily  computed (Thm.~\ref{thm: computing cpi}).   The decomposition of the Fig.~1 graph is in Fig.~2 where the $\mathcal G^6_{A, cpi}$ and $\mathcal G^6_{A, cyclic}$ entries represent, respectively, $\mathbf d^6_{A, st}$ and $\mathbf d^6_{A, cyclic}.$

 The computation of $\mathcal G^6_{A, cpi}$ follows from Def.~\ref{def: S} and Thm.~\ref{thm: computing cpi}. For instance, $S_A(A_1) = [-3+3+2-2+3]=3$ while $S_A(A_2)= [3+8+0+11+11]= 33$, so  $d_{1, 2}= \frac16 [S_A(A_1) -S_A(A_2)] = -5$ is the $\mathbf d^6_{A, st}$ entry, which is the length of the  $A_1\to A_2$ leg in $\mathcal G^6_{cpi}.$  
  The cpi subscript of $\mathcal G^6_{A, cpi}$,  which means ``closed path independent,"  is described next.

 \begin{tikzpicture}[xscale=0.5, yscale=0.52]

\draw[<-] (1.18, 0) -- (3.9, 0); 
\node[below] at (1, 0) {$A_1$};
\node[below] at (2.5, 0) {3};
\draw[->] (1.15, .06) -- (4.8, 1.9);   
\node at (1.9, .5) {\small 3};
\draw[->] (1.15, .15) -- (3.85, 3.85); 
\node at (2.1, 1.5)  {\small 2};
\draw[<-] (1, .2) -- (1, 3.9);  
\node at  (1, 3) {\small 2}; 
\node[left] at (0, 2) {$A_6$}; 
\draw[->] (.9, .1) -- (.1, 2);  
\node[left] at (.5, 1) {3};
\draw[->] (4.1, .1) -- (4.9, 1.9); 
\node[right] at  (4.5, 1) {11}; 
\node[below] at (4, 0) {$A_2$};
\draw[->] (4, .1) -- (4, 3.85);  
\node at (4, 1) {\small 11};
\draw[dashed]  (3.9, .1) -- (1.1, 3.9); 
\node at (2.1, 2.5) {\small 0};
\draw[->] (3.82, .07) -- (.23, 1.95);  
\node at (.56, 1.7) {\small 8};
\draw[<-] (4.9, 2.15) -- (4.1, 3.9);  
\node[right] at (4.5, 3) {2}; 
\node[right] at (5, 2){$A_3$};
\draw[<-] (4.8, 2.1) -- (1.24, 3.86); 
\node at (4.3, 2.35) {\small 14};
\draw[<-] (4.7, 2) -- (.2, 2); 
\node at (3.3, 2) {\small 3};
\draw[<-] (3.8, 4) -- (1.2, 4); 
\node[above] at (2.5, 4) {9};
\node[above] at (4, 4) {$A_4$}; 
\draw[<-] (3.7, 3.8) -- (.2, 2.1);  
\node at (3.05, 3.4)   {\small 1};
\node[above] at (1, 4) {$A_5$}; 
\draw[->]  (.83, 3.9) -- (.1, 2.1); 
\node[left] at (.5, 3) {8};

\node[below] at (2, -1) {{\bf a.} $\mathcal G_A^6 \approx \mathbf d_A^6$};


\draw[<-] (11.18, 0) -- (13.9, 0); 
\node[below] at (11, 0) {$A_1$};
\node[below] at (12.5, 0) {5};
\draw[->] (11.15, .06) -- (14.9, 2);   
\node at (11.9, .5) {\small 6};
\draw[->] (11.15, .15) -- (13.9, 3.9); 
\node at (12.1, 1.5)  {\small 4};
\draw[<-] (11, .2) -- (11, 3.9);  
\node at  (11, 3) {\small 5}; 
\node[left] at (10, 2) {$A_6$}; 
\draw[->] (10.9, .1) -- (10.1, 2);  
\node[left] at (10.5, 1) {3};
\draw[->] (14.1, .1) -- (14.9, 1.9); 
\node[right] at  (14.5, 1) {11}; 
\node[below] at (14, 0) {$A_2$};
\draw[->] (14, .1) -- (14, 3.85);  
\node at (14, 1) {\small 9};
\draw[dashed]  (13.9, .1) -- (11.1, 3.9); 
\node at (12.1, 2.5) {\small 0};
\draw[->] (13.82, .07) -- (10.23, 1.95);  
\node at (10.56, 1.7) {\small 8};
\draw[<-] (14.85, 2.2) -- (14.1, 3.9);  
\node[right] at (14.5, 3) {2}; 
\node[right] at (15, 2){$A_3$};
\draw[<-] (14.8, 2.1) -- (11.2, 3.86); 
\node at (14.3, 2.35) {\small 11};
\draw[<-] (14.7, 2) -- (10.25, 2); 
\node at (13.3, 2) {\small 3};
\draw[<-] (13.8, 4) -- (11.2, 4); 
\node[above] at (12.5, 4) {9};
\node[above] at (14, 4) {$A_4$}; 
\draw[<-] (13.75, 3.85) -- (10.2, 2.1);  
\node at (13.05, 3.45)   {\small 1};
\node[above] at (11, 4) {$A_5$}; 
\draw[->]  (10.83, 3.9) -- (10.11, 2.1); 
\node[left] at (10.5, 3) {8};

\node at (7.5, 2) {=};

\node[below] at (12, -1) {{\bf b.} $\mathcal G_{A, cpi}^6 \approx \mathbf d^6_{A, st}$};


\draw[->] (21.18, 0) -- (23.85, 0); 
\node[below] at (21, 0) {$A_1$};
\node[below] at (22.5, 0) {2};
\draw[<-] (21.25, .1) -- (24.85, 1.9);   
\node at (21.9, .5) {\small 3};
\draw[<-] (21.2, .2) -- (23.9, 3.9); 
\node at (22.1, 1.5)  {\small 2};
\draw[->] (21, .2) -- (21, 3.8);  
\node at  (21, 3) {\small 3}; 
\node[left] at (20, 2) {$A_6$}; 
\draw[dashed] (20.9, .1) -- (20.1, 2);  
\node[left] at (20.5, 1) {0};
\draw[dashed] (24.1, .1) -- (24.9, 1.9); 
\node[right] at  (24.5, 1) {0}; 
\node[below] at (24, 0) {$A_2$};
\draw[->] (24, .1) -- (24, 3.85);  
\node at (24, 1) {\small 2};
\draw[dashed]  (23.9, .1) -- (21.1, 3.9); 
\node at (22.1, 2.5) {\small 0};
\draw[dashed] (23.82, .07) -- (20.23, 1.95);  
\node at (20.56, 1.7) {\small 0};
\draw[dashed] (24.9, 2.1) -- (24.1, 3.9);  
\node[right] at (24.5, 3) {0}; 
\node[right] at (25, 2){$A_3$};
\draw[<-] (24.8, 2.1) -- (21.24, 3.86); 
\node at (24.3, 2.35) {\small 3};
\draw[dashed] (24.8, 2) -- (20.2, 2); 
\node at (23.3, 2) {\small 0};
\draw[dashed]  (23.8, 4) -- (21.2, 4); 
\node[above] at (22.5, 4) {0};
\node[above] at (24, 4) {$A_4$}; 
\draw[dashed] (23.75, 3.85) -- (20.2, 2.1);  
\node at (23.05, 3.4)   {\small 0};
\node[above] at (21, 4) {$A_5$}; 
\draw[dashed]  (20.83, 3.9) -- (20.11, 2.1); 
\node[left] at (20.5, 3) {0};

\node at (17.5, 2) {+};

\node[below] at (22, -1) {{\bf c.} $\mathcal G_{A, cyclic}^6 \approx \mathbf d^6_{A, cyclic}.$};

 \node[below] at (12.5, -2) {{\bf Figure 2.} Decomposition of a $\mathcal G_A^6$};

\end{tikzpicture}

Thanks to the linearity of the decomposition, the length of any closed path (not just Hamiltonian circuits) in $\mathcal G^6_A$  equals the sum of the path's lengths in $\mathcal G^6_{A, cpi}$ and $\mathcal G^6_{A, cyclic}$.  As an example, 
  the $-4$  length of  the $\mathcal G_A^6$  closed path $A_1 \stackrel{-3}{\longrightarrow}A_2 \stackrel{8}{\longrightarrow} A_6   \stackrel{3}{\longrightarrow}  A_3 \stackrel{-14}{\longrightarrow}  A_5  \stackrel{2}{\longrightarrow}  A_1$   equals the sum of its path's lengths of $A_1 \stackrel{-5}{\longrightarrow}A_2 \stackrel{8}{\longrightarrow} A_6   \stackrel{3}{\longrightarrow}  A_3 \stackrel{-11}{\longrightarrow}  A_5  \stackrel{5}{\longrightarrow}  A_1$ in $\mathcal G^6_{A, cpi}$  and  $A_1 \stackrel{2}{\longrightarrow}A_2 \stackrel{0}{\longrightarrow} A_6   \stackrel{0}{\longrightarrow}  A_3 \stackrel{-3}{\longrightarrow}  A_5  \stackrel{-3}{\longrightarrow}  A_1$ in $\mathcal G^6_{A, cyclic}$.   As this  path's length  
in  $\mathcal G^6_{A, cpi}$ is zero,    its $-4$ length  in $\mathcal G^6_{A, cyclic}$  equals the $\mathcal G^6_A$ path length.   That this closed path in $\mathcal G^6_{A, cpi}$ has length zero   is not a lucky  coincidence.  Instead, {\em all closed paths} in $\mathcal G^n_{A, cpi}$ have length zero, so these paths are  ``closed path  independent"  when computing  
  $\mathcal G^6_A$ path lengths.

\begin{Thm}\label{thm: mainA} For  $n\ge 3$,  a closed path in $\mathcal G^n_{A,  cpi}$ has length zero.  The length of a closed path in  $\mathcal G_A^n$ equals its length in $\mathcal G^n_{A, cyclic}$.   

The length of any $\mathcal G_A^n$  path  starting at $V_i$ and ending at $V_j$ equals its length in  $\mathcal G^n_{A, cyclic}$ plus the $V_i \to V_j$ length in $\mathcal G^n_{A, cpi}$. \end{Thm} 

The proof  of the first assertion 
follows immediately from the strong transitivity of the $\mathcal G^n_{A, cpi}$ components.  To check that $\mathcal G^6_{A, cpi}$ has this property, select any triplet from Fig.~2b---perhaps $\{A_1, A_5, A_3\}$; the goal is to show that $d_{1, 5}+d_{5, 3} = d_{1, 3}$.  Fom $\mathcal G^6_{A, cpi}$, this requires  $d_{1, 5}+d_{5, 3} = -5 + 11$  to equal $d_{1, 3} = 6$, which it does.   

Next, a standard induction  exercise shows that if $\mathbf d_{A, st}^n$ is strongly transitive, then {\em any path} from $A_j$ to $A_k$ has the same length as the direct  $A_j\to A_k$ path.  Reversing this last arc  creates  a closed path with length zero.  
The rest of the theorem follows immediately.  

The last assertion of Thm.~\ref{thm: mainA}  requires   the ``-1" length of   $\mathcal G^6_A$ path,  $$A_1 \stackrel{-3}{\longrightarrow} A_2 \stackrel{11}{\longrightarrow} A_4   \stackrel{2}{\longrightarrow} A_3 \stackrel{-11}{\longrightarrow} A_2  \stackrel{0}{\longrightarrow} A_5,$$ which skips $A_6$ but   visits $A_2$ twice,   to  
 equal this path's length  in  $\mathcal G^6_{A, cyclic}$ plus $-5$.  The $-5$ comes from the $A_1 \stackrel{-5}{\longrightarrow}  A_5$ arc in $\mathcal G^6_{A, cpi}$ (Fig.~2b).  This path's length in $\mathcal G^6_{A, cyclic}$ is $$A_1 \stackrel{2}{\longrightarrow} A_2 \stackrel{2}{\longrightarrow} A_4   \stackrel{0}{\longrightarrow} A_3 \stackrel{0}{\longrightarrow} A_2  \stackrel{0}{\longrightarrow} A_5$$ or 4, and, as Thm.~\ref{thm: mainA} promises,  $-1=4-5$.

\subsection{Finding the longest and shortest Hamiltonian circuits}   
According to 
Thm.~\ref{thm: mainA},  the longest/shortest  Hamiltonian path  in $\mathcal G^6_A$  (Fig.~2a) can be found by ignoring $\mathcal G^6_A$ and, instead,  analyzing  
  the reduced   $\mathcal G^6_{A, cyclic}$   with its  two    three-cycles  (Fig.~3)   
 $A_1 \stackrel{3}{\longrightarrow} A_5 \stackrel{3}{\longrightarrow} A_3 \stackrel{3}{\longrightarrow} A_1$ and $A_1 \stackrel{2}{\longrightarrow} A_2 \stackrel{2}{\longrightarrow} A_4 \stackrel{2}{\longrightarrow} A_1$.  The obvious strategy is to use these cycles as fully as possible.  To avoid premature closure by returning to a previously visited vertex, use at most two arcs of each  three-cycle.

\begin{tikzpicture}[xscale=0.5, yscale=0.52]

\draw[<-] (11.18, 0) -- (13.9, 0); 
\node[below] at (11, 0) {$A_1$};
\node[below] at (12.5, 0) {3};
\draw[->] (11.15, .06) -- (14.9, 2);   
\node at (11.9, .5) {\small 3};
\draw[->] (11.15, .15) -- (13.9, 3.9); 
\node at (12.1, 1.5)  {\small 2};
\draw[<-] (11, .2) -- (11, 3.9);  
\node at  (11, 3) {\small 2}; 
\node[left] at (10, 2) {$A_6$}; 
\draw[->] (10.9, .1) -- (10.1, 1.9);  
\node[left] at (10.5, 1) {3};
\draw[->] (14.1, .1) -- (14.9, 1.9); 
\node[right] at  (14.5, 1) {11}; 
\node[below] at (14, 0) {$A_2$};
\draw[->] (14, .1) -- (14, 3.85);  
\node at (14, 1) {\small 11};
\draw[dashed]  (13.9, .1) -- (11.1, 3.9); 
\node at (12.1, 2.5) {\small 0};
\draw[->] (13.82, .07) -- (10.23, 1.95);  
\node at (10.56, 1.7) {\small 8};
\draw[<-] (14.9, 2.1) -- (14.1, 3.9);  
\node[right] at (14.5, 3) {2}; 
\node[right] at (15, 2){$A_3$};
\draw[<-] (14.7, 2.1) -- (11.24, 3.86); 
\node at (14.3, 2.35) {\small 14};
\draw[<-] (14.7, 2) -- (10.2, 2); 
\node at (13.3, 2) {\small 3};
\draw[<-] (13.8, 4) -- (11.2, 4); 
\node[above] at (12.5, 4) {9};
\node[above] at (14, 4) {$A_4$}; 
\draw[<-] (13.75, 3.85) -- (10.2, 2.1);  
\node at (13.05, 3.4)   {\small 1};
\node[above] at (11, 4) {$A_5$}; 
\draw[->]  (10.83, 3.9) -- (10.1, 2.1); 
\node[left] at (10.5, 3) {8};

\draw[->] (21.18, 0) -- (23.85, 0); 
\node[below] at (21, 0) {$A_1$};
\node[below] at (22.5, 0) {2};
\node at (22.3, 0) {$||$};
\draw[<-] (21.25, .1) -- (24.85, 1.9);   
\node at (21.9, .5) {\small 3};
\draw[thick, <-] (21.2, .2) -- (23.9, 3.9); 
\node at (22.1, 1.5)  {\small 2};
\draw[thick, ->] (21, .2) -- (21, 3.8);  
\node at  (21, 3) {\small 3}; 
\node at (23, 1) {$||$};
\node[left] at (20, 2) {$A_6$}; 
\draw[dashed] (20.9, .1) -- (20.1, 2);  
\node[left] at (20.5, 1) {0};
\draw[dashed] (24.1, .1) -- (24.9, 1.9); 
\node[right] at  (24.5, 1) {0}; 
\node[below] at (24, 0) {$A_2$};
\draw[thick, ->] (24, .1) -- (24, 3.85);  
\node at (24, 1) {\small 2};
\draw[dashed]  (23.9, .1) -- (21.1, 3.9); 
\node at (22.1, 2.5) {\small 0};
\draw[dashed] (23.82, .07) -- (20.23, 1.95);  
\node at (20.56, 1.7) {\small 0};
\draw[dashed] (24.9, 2.1) -- (24.1, 3.9);  
\node[right] at (24.5, 3) {0}; 
\node[right] at (25, 2){$A_3$};
\draw[thick, <-] (24.8, 2.1) -- (21.24, 3.86); 
\node at (24.3, 2.35) {\small 3};
\draw[dashed] (24.8, 2) -- (20.2, 2); 
\node at (23.3, 2) {\small 0};
\draw[dashed]  (23.8, 4) -- (21.2, 4); 
\node[above] at (22.5, 4) {0};
\node[above] at (24, 4) {$A_4$}; 
\draw[dashed] (23.75, 3.85) -- (20.2, 2.1);  
\node at (23.05, 3.4)   {\small 0};
\node[above] at (21, 4) {$A_5$}; 
\draw[dashed]  (20.83, 3.9) -- (20.11, 2.1); 
\node[left] at (20.5, 3) {0};

\node at (17.5, 2) {$\to$};

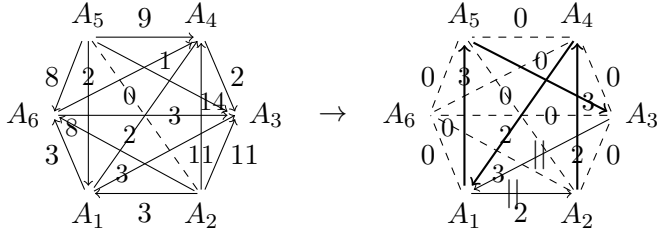
\node[below] at (17.4, -1.2) {{\bf Figure 3.} Finding $\mathcal G_A^6$'s longest Hamiltonian circuit via $\mathcal G^6_{A, cyclic}$.};

\end{tikzpicture}

Consequently,  starting at $A_1$, use  only the edges $A_1 \stackrel{3}{\longrightarrow} A_5 \stackrel{3}{\longrightarrow} A_3 $  from the first three-cycle; this  is why  $A_3  \stackrel{3}{\longrightarrow} A_1$ arc  is crossed out in Fig.~3.  As only one arc can leave $A_1$,   this   eliminates  the crossed-out  $A_1 \stackrel{2}{\longrightarrow} A_2$ arc.   
Thus,  the Hamiltonian circuit depends upon the four bold $\mathcal G^6_{A, cyclic}$ arcs in Fig.~3.  Finding the longest Hamiltonian circuit now is immediate---use zero-length arcs to ensure all vertices are visited and to connect the two basic ones.  An answer is the Hamiltonian path 
$$ (A_1 \stackrel{3}{\longrightarrow} A_5 \stackrel{3}{\longrightarrow} A_3)   \stackrel{0}{\longrightarrow} A_6  \stackrel{0}{\longrightarrow} (A_2  \stackrel{2}{\longrightarrow} A_4  \stackrel{2}{\longrightarrow} A_1)$$ of length 10.  The reversal of this closed path has length $-10$.  It is clear from the graph's structure that  these are the longest and shortest $\mathcal G^6_{A, cyclic}$ Hamiltonian circuits, so (Thm.~3)  they define the longest and shortest  $\mathcal G^6_A$ Hamiltonian circuits with the same lengths.

To appreciate what is going on, notice that what complicates  computing the Fig.~2a  path length $A_5 \stackrel{14}{\longrightarrow} A_3 \stackrel{-3}{\longrightarrow} A_1 \stackrel{-2}{\longrightarrow} A_5$ of $14-3-2=9$ are the subtractions/cancelations.  To see what they are, let $u$, $v$, and $w$ be the  canceled terms, respectively, for arcs $\widehat{A_5A_3}, \widehat{A_3A_1}$, and $\widehat{A_1A_5}$.  That is, after cancellations,  the path length computation would be $(14-u)+(-3-v) +(-2-w)= 9$ where  $u+v+w=0$,  which satisfies cpi from Eq.~\ref{def: st}.  That is, the cancelations in computing closed path lengths are linear expressions that define a   $\mathbb {ST}_A^n$ graph.  The optimal choice of cancelations for a $\mathcal G^n_A$ is  the $\mathbb {ST}_A^n$ graph that most closely approximates  $\mathcal G^n_A$, which  is  
  its orthogonal projection into  $\mathbb {ST}_A^n$, or $\mathcal G^n_{A, cpi}$.  Indeed, with Fig.~2b, the optimal cancelation is from $\mathcal G^6_{A, cpi}$ where $u+v+w= 11-6-5=0.$  Removing these cancelations leaves $\mathcal G^6_{A, cyclic}$ with the path length $A_5 \stackrel{3}{\longrightarrow} A_3 \stackrel{3}{\longrightarrow} A_1 \stackrel{3}{\longrightarrow} A_5$ of length $9$ where no modifications are needed in the computations.

\subsection{The symmetric case}  For the standard  symmetric TSP, the distance between vertices is the same in each direction.  Thus  the $d_{i, j} = -d_{j, i}$ asymmetric condition is replaced with  the $d_{i, j} = d_{j, i}$ requirement.  This  requires a different decomposition of the new $\mathbb R_S^{n\choose2}$ where, with the orthogonal projection, leads to a different type of  inherent symmetries.  Otherwise, everything remains much the same.   Of importance, the decomposition  provides a best possible simplification   for closed paths.   Details for the general symmetric and asymmetric cases  are given elsewhere \cite{new}, but a flavor of what happens is described next.

A $n=6$ example is given in Fig.~4; recall,  the S subscript  in $\mathcal G^6_S$ refers to ``symmetric."  While a typical goal is to find the shortest Hamiltonian circuit, everything  extends to the analysis of any path, whether  closed or not.  


\begin{tikzpicture}[xscale=0.5, yscale=0.6]

\draw[<-] (1.18, 0) -- (3.9, 0); 
\node[below] at (1, 0) {$A_1$};
\node[below] at (2.5, 0) {1};
\draw (1.15, .06) -- (4.9, 2);   
\node at (1.9, .5) {\small 5};
\draw (1.15, .15) -- (3.9, 3.9); 
\node at (2.1, 1.5)  {\small 3};
\draw (1, .2) -- (1, 3.9);  
\node at  (1, 3) {\small 4}; 
\node[left] at (0, 2) {$A_6$}; 
\draw (.9, .1) -- (.1, 2);  
\node[left] at (.5, 1) {1};
\draw (4.1, .1) -- (4.9, 1.9); 
\node[right] at  (4.5, 1) {7}; 
\node[below] at (4, 0) {$A_2$};
\draw (4, .1) -- (4, 3.85);  
\node at (4, 1) {\small 7};
\draw  (3.9, .1) -- (1.1, 3.9); 
\node at (2.1, 2.5) {\small 3};
\draw (3.82, .07) -- (.23, 1.95);  
\node at (.56, 1.7) {\small 4};
\draw (4.9, 2.1) -- (4.1, 3.9);  
\node[right] at (4.5, 3) {2}; 
\node[right] at (5, 2){$A_3$};
\draw (4.9, 2.05) -- (1.24, 3.86); 
\node at (4.3, 2.35) {\small 3};
\draw (4.8, 2) -- (.2, 2); 
\node at (3.3, 2) {\small 5};
\draw (3.8, 4) -- (1.2, 4); 
\node[above] at (2.5, 4) {1};
\node[above] at (4, 4) {$A_4$}; 
\draw (3.75, 3.85) -- (.2, 2.1);  
\node at (3.05, 3.4)   {\small 7};
\node[above] at (1, 4) {$A_5$}; 
\draw (.83, 3.9) -- (.1, 2.1); 
\node[left] at (.5, 3) {7};

\node[below] at (2, -1) {{\bf Figure 4.}  A symmetric $\mathcal G_S^6 \in \mathbb G^6_S$};

\end{tikzpicture}

Central to the discussion is the average Hamiltonian path length in $\mathcal G^n_S$ denoted by  $T(\mathcal G^n_S)$.  

\begin{Def}\label{def: symmetric}  For  a graph $\mathcal G^n_S$, let $S_S(A_j)$ be the sum of arc lengths attached to vertex  $A_j$, $j=1, \dots, n.$  Let $T(\mathcal G^n_S) = \frac1{n-1} \sum_{j=1}^n S_S(A_j)$.  \end{Def} 

So $S_S(A_j)$ (Def.~\ref{def: symmetric}) and $S_A(A_j)$ (Eq.~\ref{def: S}) are the same.  That  $T(\mathcal G^n_S)$ is the average Hamiltonian path length follows from the fact that $\frac1{n-1} S_S(A_j)$  is the average length of the $(n-1)$ arcs attached to $A_j$.   
Illustrating  with Fig.~4, $S_S(A_1) = 1+5+3+4+1= 14, S_S(A_2) = 22, S_S(A_3) =  22, S_S(A_4) = 20, S_S(A_5) = 18, S_S(A_6)= 24$, so $T(\mathcal G^6_S) = \frac15(120) = 24.$   
When considering non-Hamiltonian  closed paths, restrict the $S_S(A_j)$ values and the summation defining $T$ to the relevant arcs and vertices.

\begin{tikzpicture}[xscale=0.5, yscale=0.6]

\draw[<-] (1.18, 0) -- (3.9, 0); 
\node[below] at (1, 0) {$A_1$};
\node[below] at (2.5, 0) {1};
\draw (1.15, .06) -- (4.9, 2);   
\node at (1.9, .5) {\small 5};
\draw (1.15, .15) -- (3.9, 3.9); 
\node at (2.1, 1.5)  {\small 3};
\draw (1, .2) -- (1, 3.9);  
\node at  (1, 3) {\small 4}; 
\node[left] at (0, 2) {$A_6$}; 
\draw (.9, .1) -- (.1, 2);  
\node[left] at (.5, 1) {1};
\draw (4.1, .1) -- (4.9, 1.9); 
\node[right] at  (4.5, 1) {7}; 
\node[below] at (4, 0) {$A_2$};
\draw (4, .1) -- (4, 3.85);  
\node at (4, 1) {\small 7};
\draw  (3.9, .1) -- (1.1, 3.9); 
\node at (2.1, 2.5) {\small 3};
\draw (3.82, .07) -- (.23, 1.95);  
\node at (.56, 1.7) {\small 4};
\draw (4.9, 2.1) -- (4.1, 3.9);  
\node[right] at (4.5, 3) {2}; 
\node[right] at (5, 2){$A_3$};
\draw (4.9, 2.05) -- (1.24, 3.86); 
\node at (4.3, 2.35) {\small 3};
\draw (4.8, 2) -- (.2, 2); 
\node at (3.3, 2) {\small 5};
\draw (3.8, 4) -- (1.2, 4); 
\node[above] at (2.5, 4) {1};
\node[above] at (4, 4) {$A_4$}; 
\draw (3.75, 3.85) -- (.2, 2.1);  
\node at (3.05, 3.4)   {\small 7};
\node[above] at (1, 4) {$A_5$}; 
\draw (.83, 3.9) -- (.1, 2.1); 
\node[left] at (.5, 3) {7};


\node at (7.4, 2) {$\to$};

\draw (11.18, 0) -- (13.9, 0); 
\node[below] at (11, 0) {$A_1$};
\node[below] at (12.5, 0) {-2};
\draw (11.15, .06) -- (14.9, 2);   
\node at (11.9, .5) {\small 2};
\draw (11.15, .15) -- (13.9, 3.9); 
\node at (12.1, 1.5)  {\small 0.5};
\draw (11, .2) -- (11, 3.9);  
\node at  (11, 3) {\small 2}; 
\node[left] at (10, 2) {$A_6$}; 
\draw (10.9, .1) -- (10.1, 2);  
\node[left] at (10.5, 1) {-2.5};
\draw (14.1, .1) -- (14.9, 1.9); 
\node[right] at  (14.5, 1) {2}; 
\node[below] at (14, 0) {$A_2$};
\draw (14, .1) -- (14, 3.85);  
\node at (13.9, 1) {\small 2.5};
\draw  (13.9, .1) -- (11.1, 3.9); 
\node at (12.1, 2.5) {\small -1};
\draw (13.82, .07) -- (10.23, 1.95);  
\node at (10.56, 1.7) {\small -1.5};
\draw (14.9, 2.1) -- (14.1, 3.9);  
\node[right] at (14.5, 3) {-2.5}; 
\node[right] at (15, 2){$A_3$};
\draw[<-] (14.9, 2.05) -- (11.24, 3.86); 
\node at (14.3, 2.35) {\small -1};
\draw  (14.8, 2) -- (10.2, 2); 
\node at (13.3, 2) {\small -0.5};
\draw    (13.8, 4) -- (11.2, 4); 
\node[above] at (12.5, 4) {-2.5};
\node[above] at (14, 4) {$A_4$}; 
\draw  (13.75, 3.85) -- (10.2, 2.1);  
\node at (13.05, 3.4)   {\small 2};
\node[above] at (11, 4) {$A_5$}; 
\draw  (10.83, 3.9) -- (10.11, 2.1); 
\node[left] at (10.5, 3) {2.5};

\node at (17.5, 2) {+};

\node at (19.5, 2) {$ T(\mathcal G^6_S) $};

\node[below] at (12, -1) { $\mathcal G_{S, cyclic}^6 \in \mathbb C^6_S$};
\node[below] at (2, -1) {$\mathcal G^6_{S, cyclic}$};

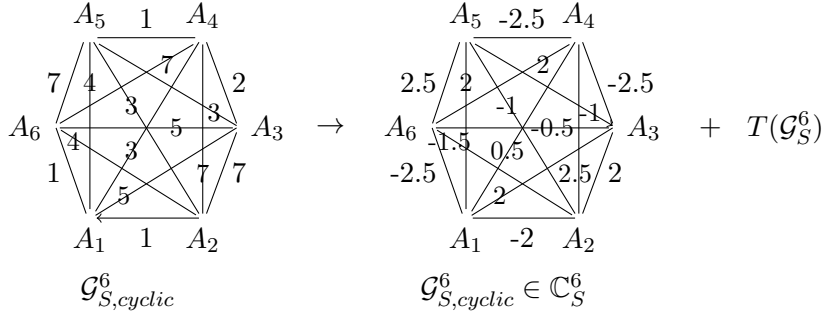
 \node[below] at (8.5, -2) {{\bf Figure 5.} Decomposition of a $\mathcal G_S^6 \in \mathbb G^6_S$};

\end{tikzpicture}

The scheme replaces entries of $\mathcal G^n_S$ with values that, after removing irrelevant terms (which are  similar to $\mathcal G^n_{A, cpi}$),  can be viewed as differences from the average arc length, as in Fig.~5.  Thus negative values represent ``smaller than average" costs.  The theorem is that {\em the  length of a Hamiltonian circuit in $\mathcal G^n_S$  is the length of its path in $\mathcal G^n_{S, cycle}$ plus $T(\mathcal G^n_S)$.}  Standard modifications handle paths that are not closed and/or incomplete graphs.

Each  $\mathcal G^6_{S, cyclic}$ vertex has an arc with negative (i.e., below average) length; this assertion  holds  more generally   for any $\mathcal G^n_{S, cyclic}$.  By observation and using arcs with negative values as often as possible, an Hamiltonian  path in $\mathcal G^6_{S, cyclic}$ with shortest length of -11  is 
\begin{equation}\label{tour} A_1 \stackrel{-2.5}{\longrightarrow} A_6 \stackrel{-0.5}{\longrightarrow}A_3 \stackrel{-2.5}{\longrightarrow} A_4 \stackrel{-2.5}{\longrightarrow} A_5 \stackrel{-1}{\longrightarrow} A_2 \stackrel{-2}{\longrightarrow} A_1.\end{equation}  This  defines the  
  shortest Hamiltonian path in $\mathcal G^6_S$ with length $24-11=13.$

In general,  the shortest  $\mathcal G^n_S$ Hamiltonian path  is bounded below by $T(\mathcal G^n_S)$ plus the sum of the  $n$ smallest  $\mathcal G^n_{S, cyclic}$ arc lengths.  With $\mathcal G^6_S$, this is $24+ (-2.5 -2.5 -2.5 -2 - 1.5 -1) = 12.$ This shortest Hamiltonian path  length is larger than,  but close to, its lower bound.  

 Although the   $\mathcal G^6_S$ entries can be random numbers,  the advantage of the  $\mathcal G^6_{S, cyclic}$ entries is that each entry  has a meaning with respect to  a path's length.  Thus,  a  way to find the shortest Hamiltonian path is to  rank the arcs according to length, where ``smaller is better."  Start with the first n arcs to determine whether they form a path. If not, then iteratively add  arcs to see whether it creates a Hamiltonian path.   
 To illustrate, the $\mathcal G^6_{S, cyclic}$   arcs with negative values are

\begin{equation} \begin{array} {c|ccc|c|cc} \textbf{Length} & \textbf{Arc} &  \textbf{Arc} &  \textbf{Arc} &  \textbf{Length} & \textbf{Arc} &  \textbf{Arc}\\ \hline 
-2.5 & \widehat{A_1A_6}, &
 \widehat{A_3A_4}, &
\widehat{A_4A_5} &
-2 & \widehat{A_1A_2} \\   
-1.5 & \widehat{A_2A_6} & &  &
-1 & \widehat{A_3A_5}, & \widehat{A_2A_5} \\  -0.5 & \widehat{A_3A_6}

\end{array}\end{equation}  This array emphasizes using $\widehat{A_1A_6}, \,  \widehat{A_3A_4}, \widehat{A_4A_5}, $ where 
$\widehat{A_3A_6} $ forms a transition between the first and the second arc, and $\widehat{A_5A_2}$ includes the missing $A_2$, while $\widehat{A_2A_1}$ completes the journey.  This is the Eq.~\ref{tour} tour.  \smallskip


\begin{tikzpicture}[xscale=0.5, yscale=0.5]

\draw[thick] (-1, 0) -- (5, 0);
\node[right] at (5, .1) {$\mathbb{ST}^n_A$; Voting};

\draw[thick] (0, -1) -- (0, 5); 
\node[above] at (0, 5.0) {$\mathbb C^n_A$: TSP and graphs}; 

\node at (4, 4) {$\bullet$};
\node[right] at (4, 4) {$\mathbf d_A^n \in \mathbb R_A^{n\choose2}$};
\draw[->] (4, 3.9) -- (4, 0.1);
\draw[->] (3.9, 4) -- (0.1, 4); 
\node[below] at (4, -.1) {Borda outcome};
\node[left] at (0, 4)  {$\mathcal G^n_{A, cyclic}$};

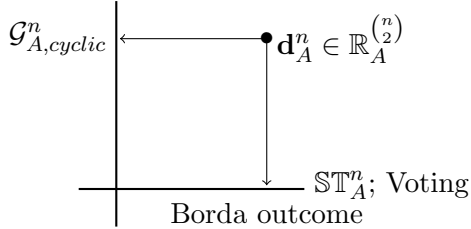
\node at (2.5, -2) {{\bf Figure 6.} Resolutions via projections};

\end{tikzpicture}

\section{Conclusion} Figure 6 summarizes much of what was discussed.  Namely, both the voting problem and  TSP can be described in terms of  a point $\mathbf d_A^n \in \mathbb R^{n\choose2}_A$.   This space has an orthogonal decomposition into  the strongly transitive linear subspace $\mathbb {ST}^n_A$ and the  linear subspace $\mathbb C^n_A$ consisting of three-cycles. The $\mathbf d_{A, st}^n$  component of  $\mathbf d_A^n$ in the $\mathbb {ST}^n_A$ direction  resolves several concerns from voting theory, but frustrates the analysis of paths in TSP.  Conversely, $\mathbf d_A^n$'s component in the $\mathbb C^n_A$ direction, $\mathbf d_{A, cyclic}^n$,  determines all closed path properties for the TSP, but erects barriers and paradoxes for voting theory.  A resolution   is to project  $\mathbf d_A^n$  into the appropriate subspace that has  positive properties for the problem being considered.


  To introduce  the rest of the material in this concluding section, let me brag about my two granddaughters.  Heili has a deep interest in neurobiology and gymnastics, while Tatjana is heavily involved in the social sciences and ballet.  Who  is  older?  
It is impossible to tell from this information;  comparing ages requires having  relevant inputs.  More generally,  to realize any specified objective,   for the data to be useful; they  must, in some manner,  be directed toward the stated purpose.  This truism extends to voting theory and the TSP.   

To see this with Arrow's result, return  to the introductory voting example where the alternatives are three cities.  Suppose the 25 voters with  the $A_1\succ A_2 \succ A_3$ ranking judge cities according the ease of finding parking,  the 23 who prefer $A_2\succ A_3 \succ A_1$ evaluate cities according to the popularity  of local professional  basketball teams, and the  22 preferring $A_3\succ A_1\succ A_2$  assess them according to  available types of craft beer.  Thus the $A_1\succ A_2$ outcome  by 47:23 is a conglomeration of  attitudes about parking,  basketball, and beer. Of importance,  {\em nothing} is directed toward Arrow's explicitly stated  objective of having a transitive ranking.  This means (with these apple and orange comparisons)  that rather than a surprise,   non-transitive outcomes should   be  anticipated.  Indeed, expect a transitive outcome  {\em only if}  the  
 inputs  contribute toward this transitivity target.  
This comment leads to the decomposition (Fig.~6);  the $\mathbf d^n_{A, cyclic}$ component of $\mathbf d_A^n$ that runs counter to the stated objective is identified and  dismissed. This leads to the approach described in the first paragraph of the proof of Thm.~\ref{thm: computing cpi}, which is called the ``summation method" in \cite{geometry}.   

Another way to handle this problem is to express the inputs, the preferences of voters, in terms of the transitivity   objective.  To do so,  describe  the  $A_1\succ A_2 \succ A_3$ ranking  as $(A_1\succ A_2, 1), (A_1\succ A_3, 2)$ and $(A_2\succ A_3, 1)$, which  restates this    ranking in terms of  pairs that now  are   in  a strongly transitive format.  Here $d_{1, 2}=1$ captures that $A_1$ is ranked one spot above $A_2$, while $d_{1, 3}=2$ means that $A_1$ is ranked two spots above $A_3$.  Choosing the triplet ($A_1, A_2, A_3$) requires checking  whether $d_{1, 2} + d_{2, 3} = 1+1$ equals $d_{1, 3}=2,$ which it does. By using this strongly transitive approach, which is called ``Intensity of IIA" or IIIA (introduced in \cite{geometry}), the profile information for  $\{A_1, A_3\}$ has 25 voters with $(A_1\succ A_3, 2)$, 23 with $(A_1\succ A_3, -1)$ and 22 with $(A_1\succ A_3, -1)$,  The IIIA tally for  $\{A_1, A_3\}$   sums the products of the number of voters with each ranking times its intensity.  Here the   tally is $25(2) + 23(-1) + 22(-1) = 5$ leading to  $A_1\succ A_3$, rather than the above $A_3\succ A_1$ (Eq.~\ref{eq: introductory}) that forced a cycle.  Even stronger, the IIIA tallies over the three pairs not only are transitive, they are strongly transitive.  
 
 With  any number of alternatives, the IIIA outcome for the majority vote always is strongly transitive.  It must be because $\mathbb{ST}_A^n$ is a linear subspace, and the outcome is a summation of   strongly transitive profiles.   This IIIA  method is equivalent to the Borda Count \cite{geometry}; a conclusion   that should be anticipated (particularly with Thm.~\ref{thm: BC}).   
  As shown in \cite{geometry},  by replacing IIA with IIIA, Arrow's Theorem now has a positive conclusion.  That is,  
  by using  data that is consistent with the objective of transitivity,  the problems of Arrow's Theorem disappear.   

Similarly with TSP, each difference from the average cost between two vertices, as catalogued with $\mathcal G^6_A$ (Fig.~1),  could be based on different attributes; e.g., the $d_{1, 2}$  difference between $A_1$ and $A_2$ might reflect the  topography while $d_{2, 3}$  between $A_2$ and $A_3$ could be caused by traffic restrictions.  As the objective concerns path lengths of closed curves, the goal must be to use inputs that contribute to the specified goal.  Here the strongly transitive component $\mathbf d_{A, cpi}^n$ pushes the outcome in a linear fashion that conflicts with the goal of  having circular paths.  Dismissing this $\mathbf d_{A, cpi}^n$  term (as indicated in Fig.~6)  leaves $\mathbf d^n_{A, cyclic}$ data, which are consistent with the   global objective of finding closed path properties.

\section{Proofs}\label{sect: proof} 

Beyond proving Thm.~3, an  intent of this section is to demonstrate  the above tools.\smallskip

\noindent{\em Proof of Thm.~3.}  Assume that the paired rankings defined by  $\mathbf d^n$  have $A_1$ and $A_n$, respectively, the unique Condorcet winner and  
 loser.  As $A_1$  is the unique Condorcet winner,  the associated $\mathbf d_A^n$  must have  $d_{1, k}>0$ for $k=2, \dots, n$.  Similarly, as $A_n$ is the unique Condorcet loser, it must be that  $d_{n, k}<0$ for all $k=1, \dots, n-1.$   
  From Eq.~\ref{def: S}, this means that $S_A(A_1)>0$ and $S_A(A_n)<0$.  According to Eq.~\ref{eq: StoB}, $A_1$ is Borda ranked above $A_n$.  As the Borda ranking is the  $\mathbf d^n_{A, st},$  ranking, this proves the  theorem. Notice,     the Condorcet uniqueness conditions are unnecessary  conveniences.   $\square$  \smallskip  

It remains to show for $n=3, 4$  that if $A_1$ and $A_n$ are, respectively, the top and bottom ranked alternative for $\mathbf d^n_{A, st}$, and if $\mathbf d_A^n = \mathbf d^n_{A, st} + \mathbf d^n_{A, cyclic}$  is transitive,  then $A_1\succ A_n$.  The approach uses the fact, which follows from strong transitivity,  that $d_{1, n}$ for $\mathbf d^n_{A, st}$ is an upper bound for all other $d_{s, u}$ values.  This size of $d_{1, n}$  requires  the cyclic perturbations to affect and reverse  smaller $d_{u, s}$ values, which   
 violates transitivity before they can  impact on the $d_{1, n}$ term to reverse $A_1\succ A_n$.  Assume the $\mathbf d^n_{A, st}$ ranking is $A_1\succ A_2\succ \dots  \succ A_{n-1} \succ A_n$.\smallskip

The proof is immediate for $n=3$.   All cyclic $n=3$ vectors  are   $\alpha$ multiples  of $(1, -1; 1)$, so  $\mathbf d_A^3=\mathbf d^3_{A, st} + \mathbf d^3_{A, cyclic}=  
 (d_{1, 2} + \alpha, d_{1, 3} -\alpha; d_{2, 3} + \alpha).$   If $\alpha\ge0$,  then the  $A_1 \succ A_2$ and   $A_2 \succ A_3$ rankings remain unchanged.  The  $A_1>A_3$ ranking persists as long as $d_{1, 3} - \alpha>0$.  As soon as  $\alpha>  d_{1, 3}$,  the $A_1\succ A_3$  ranking reverses to become $A_3\succ A_1$, which converts the set of pairwise  rankings from  transitive to  cyclic.   

For all  $\alpha<0$ values, $A_1\succ A_3$.  But the system becomes cyclic as soon as   
 $-\alpha$ equals the second largest of $\{d_{1, 2}, d_{2, 3}\}.$ For instance, suppose $d_{1, 2} < d_{2, 3}.$  For $\alpha=-d_{1, 2}$, the rankings are $A_1\succ A_3, A_2\succ A_3,$ and $A_2\sim A_1$ defining  the transitive   $A_1\sim A_2\succ A_3$.  But once $\alpha=-d_{2, 3}$, the rankings are the non-transitive $A_1\succ A_3, A_3\sim A_2, A_2\succ A_1$. Stated differently, if $\mathbf d_A^3$ has a transitive ranking, then  $A_1\succ A_3$, which proves the assertion. \smallskip  

The  theme of the proof for $n\ge 4$   is that, because of strong transitivity,  $d_{1, n}$ from $\mathbf d^n_{A, st}$ is so large  that  before  the cyclic components can reverse  the $A_1\succ A_n$ ranking, they change enough rankings of other pairs to  define a cyclic outcome  that  violates the transitivity of $\mathbf d_A^n$.   One approach  follows: \begin{enumerate}
\item Assume that the cyclic terms force $A_n\succ A_1$.  It must be shown that the pairs do not define a transitive outcome.
\item  The size of the perturbations of cyclic terms that accompany  the    $A_n\succ A_1$ assumption requires some other alternative, say $A_k$, to satisfy $A_k\succ A_n$.  If $\mathbf d_A^n$ is transitive, as assumed, then  $A_k\succ A_n\succ A_1,$  and $A_k\succ A_1.$ \item
The size of the cyclic perturbations  that cause $A_k\succ A_1$ forces some other  alternative, $A_j$, to satisfy $A_j\succ A_k$.  This requires  $A_j\succ A_k\succ A_n\succ A_1$, or   $A_j\succ A_1$.
\item The next step is to  show that this $A_j\succ A_1$ ranking drives the size of certain cyclic terms to be large enough so that $A_n\succ A_j$
\item Step 4 is the sought  contradiction.  The $A_j\succ A_n$ ranking follows from the assumption of the transitivity of $\mathbf d_A^n$, and the conflicting    $A_n\succ A_j$ is a direct result of computations based on the impacts of the cyclic terms.  Thus    the assumption that $\mathbf d_A^n$ is transitive is false, which proves the conclusion.
\end{enumerate}

To illustrate   this program with  $n=4$,   all of the following $d_{i, j}$ values  come  from $\mathbf d^4_{A, st}$.   Assume that $\mathbf d_A^4$ is transitive.  Using the Thm.~2 basis for $\mathbb C_A^4$, where $\alpha_{j, k}$ is the coefficient for $\mathbf d^4_{1, j, k}$,   a representation for $\mathbf d_A^4= \mathbf d^4_{A, st} + \mathbf d^4_{A, cyclic}$ is   
  \begin{equation}\label{eq: forconclusion} \mathbf d_A^4= ({d_{1, 2} + \alpha_{2, 3} + \alpha_{2, 4},  d_{1, 3} - \alpha_{2, 3} + \alpha_{3, 4}, d_{1, 4} -\alpha_{24} - \alpha_{3, 4}; d_{2, 3} + \alpha_{2, 3},  d_{2, 4} + \alpha_{2, 4},  d_{3, 4} + \alpha_{3, 4}}).\end{equation} 

According to Eq.~\ref{eq: forconclusion},  
  if  $A_4\succ A_1$, then   $d_{1, 4} -\alpha_{2, 4} - \alpha_{3, 4} <0$, or $0<d_{1, 4} < \alpha_{2, 4} +\alpha_{3, 4}$.  This forces  one or both of  $\alpha_{3, 4}$, $\alpha_{2, 4}$ to be positive.   
 \smallskip
 
 Step 2.   Start with the assumption that $\alpha_{3, 4}>0$, which requires from Eq.~\ref{eq: forconclusion} that $A_3>A_4$.  From  the transitivity of $\mathbf d_A^4$, this  requires $A_3>A_1$, or $d_{1, 3} -\alpha_{2, 3} + \alpha_{3, 4}<0,$ which  leads to the inequality \begin{equation}\label{eq: 23} 0< d_{1, 3} + \alpha_{3, 4} <\alpha_{2, 3}.\end{equation}   
 
 Step 3.  Combining $\alpha_{2, 3}>0$  (Eq.~\ref{eq: 23}) with  $\mathbf d_A^4$'s transitivity (and Eq.~\ref{eq: forconclusion}),  requires  $A_2>A_3>A_4> A_1$.  The $A_2>A_1$ ranking forces  (Eq.~\ref{eq: forconclusion}) $\alpha_{2, 3} + \alpha_{2, 4}<-d_{1, 2}<0,$  which,  with Eq.~\ref{eq: 23},  becomes  \begin{equation}\label{eq: alpha} \alpha_{2. 4} < -d_{1, 2} - \alpha_{2, 3} <  -d_{1, 2} - d_{1, 3} - \alpha_{3, 4}<0.\end{equation} 
 
 Step 4.  As $\alpha_{2, 4}<0$ (Step 3), it follows from $A_4>A_1$ and $d_{1, 4} -\alpha_{2, 4} - \alpha_{3, 4}<0$  that \begin{equation} \label{eq: large} \alpha_{3, 4} > d_{1, 4}\ge d_{j,k}.\end{equation}  The last inequality follows because $d_{1, 4}$ is the largest $d_{j, k}$ value (strong transitivity).   
 
 Step 5.  By using Eq.~\ref{eq: alpha}, the computation for the   $\{A_2, A_4\}$ ranking depends on the sign of  $d_{2, 4} + \alpha_{2, 4} <  - d_{1, 2} -d_{1, 3} +[ -\alpha_{3, 4} + d_{2, 4}]$.  Because  $-\alpha_{3, 4} + d_{2, 4}<0$ (Eq.~\ref{eq: large}), it follows that  $d_{2, 4} + \alpha_{2, 4}<0$, which is $A_4\succ A_2$ and the desired contradiction.   
 \medskip

The remaining case of $\alpha_{2, 4}>0, \alpha_{3, 4}\le 0$ is simpler because from $A_4\succ A_1$  (Eq.~\ref{eq: forconclusion})  \begin{equation}\label{eq: 24large} \alpha_{2, 4}>d_{1, 4}\ge d_{j, k}.\end{equation}  
 
 Step 2.  According to Eq.~\ref{eq: forconclusion}, the $\alpha_{2, 4}>0$ inequality and the transitivity of $\mathbf d^4$  mandates   $A_2>A_4>A_1$, or from $A_2\succ A_1$  that $d_{1, 2} + \alpha_{2, 3} +\alpha_{2, 4} <0.$ or that $\alpha_{2, 3} < -d_{1, 2} -\alpha_{2, 4}<0.$  
 
 Step 3.  Using this inequality and Eq.~\ref{eq: 24large}, the $\{A_2, A_3\}$ ranking equation is determined by 
 $d_{2, 3} + \alpha_{2, 3} < -d_{1, 2} + [d_{2, 3}-\alpha_{2, 4}] <0$, or $A_3>A_2>A_4>A_1$.
 
 Step 4.  According to the $A_3>A_1$ ranking,   $d_{1, 3} - \alpha_{2, 3} + \alpha_{3, 4}<0 $ or (Step 2)  $d_{1, 3}  + \alpha_{3, 4}<\alpha_{2, 3} < -d_{1, 2} -\alpha_{2, 4}<0.$ What follows uses  
  the inequality $\alpha_{3, 4} < -d_{1, 3} -d_{1, 2} - \alpha_{2, 4}.$
 
Step 5.  The  
 $\{A_3, A_4\}$ ranking depends on the sign of  $d_{3, 4} + \alpha_{3, 4} < -d_{1, 2} - d_{1, 3}  +[-\alpha_{2, 4}+d_{3, 4}] <0$.  As the term in the brackets is negative, the value is negative, leading to the  $ A_4>A_3$ contradiction.  This completes the proof for $n=4$.

\end{document}